\documentclass[11pt]{article}
\usepackage{mathrsfs}
\usepackage{amssymb}
\usepackage{amsmath}
\usepackage{amsfonts}
\usepackage{graphicx}

\newtheorem{lemma}{Lemma}[section]
\newtheorem{theorem}{Theorem}[section]
\newtheorem{remark}{Remark}[section]

\date{August 29 \rm 2011}

\numberwithin{equation}{section}

\def\dint{\displaystyle\int}

\def\dbE{{\mathop{\rm l\negthinspace E}}}

\def\dbP{{\mathop{\rm l\negthinspace P}}}
\def\dbR{{\mathop{\rm l\negthinspace R}}}
\def\={\buildrel \triangle \over =}

\begin{document}

\title{\bf Mean-Field Backward Doubly Stochastic Differential Equations
and Applications\footnote{This work is supported in part by National
Natural Science Foundation of China (Grants 10771122 and 11071145),
Natural Science Foundation of Shandong Province of China (Grant
Y2006A08), Foundation for Innovative Research Groups of National
Natural Science Foundation of China (Grant 10921101), National Basic
Research Program of China (973 Program, No. 2007CB814900),
Independent Innovation Foundation of Shandong University (Grant
2010JQ010), Graduate Independent Innovation Foundation of Shandong
University (GIIFSDU).}
 }
\author{Tianxiao Wang\footnote{School of Mathematics, Shandong
University, Jinan 250100, China},~~Qingfeng Zhu\footnote{ 1.School
of Mathematics, Shandong University, Jinan 250100, China; 2.School
of Statistics and Mathematics, Shandong University of Finance and
Economics, Jinan 250014, China},~~and~~Yufeng
Shi\footnote{Corresponding author. Email:yfshi@sdu.edu.cn. Shandong
University, Jinan 250100, China}
 }

 \maketitle

\begin{abstract}

 Mean-field backward doubly stochastic differential equations
(MF-BDSDEs, for short) are introduced and studied. The existence and
uniqueness of solutions for MF-BDSDEs is established. One
probabilistic interpretation for the solutions to a class of
nonlocal stochastic partial differential equations (SPDEs, for
short) is given. A Pontryagin's type maximum principle is
established for optimal control problem of MF-BDSDEs. Finally, one
backward linear quadratic problem of mean-field type is discussed to
illustrate the direct application of above maximum principle.

\end{abstract}

\vspace{0.3cm} \bf Keywords. \rm Mean-field backward doubly
stochastic differential equations, nonlocal stochastic partial
differential equations, maximum principle.

\vspace{0.3cm}


\section{Introduction}

Backward doubly stochastic differential equation (BDSDE for short)
of the form
\begin{eqnarray*}
Y_t=\xi +\dint_t^Tf(s,Y_s,Z_s)ds+\dint_t^Tg(s,Y_s,Z_s) d
\overleftarrow B_s-\dint_t^TZ_s d \overrightarrow W_s ,\ 0\leq t\leq
T,
\end{eqnarray*}
was firstly initiated by Pardoux--Peng \cite{Pardoux-Peng 1994} to
give probabilistic interpretation for the solutions of a class of
semilinear stochastic partial differential equations (SPDEs for
short). BDSDEs have not only emerged as a natural and convenient
tool in the context of SPDEs, see Bally--Matoussi
\cite{Bally-Matoussi 2001}, Hu--Ren \cite{Hu-Ren 2009},
Pardoux--Peng \cite{Pardoux-Peng 1994}, Ren--Lin--Hu
\cite{Ren-Lin-Hu 2009}, Zhang--Zhao (\cite{Zhang-Zhao 2007},
\cite{Zhang-Zhao 2010}), but also recently gained interest in other
fields as well, especially in relations to the stochastic optimal
control problems, see Bahlali--Gherbal \cite{Bahlali-Gherbal 2010},
Han--Peng--Wu \cite{Han-Peng-Wu 2010}, Zhang--Shi \cite{Zhang-Shi
2010}.

{\it McKean--Vlasov} stochastic differential equation of the form
\begin{eqnarray}\label{eq:1.1}
dX(t)=b(X(t),\mu(t))dt+dW(t),\quad t\in[0,T],\quad X(0)=x,
\end{eqnarray}
where $$b(X(t),\mu(t))=\int_\Omega
b(X(t,\omega),X(t;\omega'))\dbP(d\omega')=\dbE[b(\xi,X(t))]
\Big|_{\xi=X(t)},$$ $b:\dbR^m\times\dbR\to\dbR$ being a (locally)
bounded Borel measurable function and $\mu(t;\cdot)$ being the
probability distribution of the unknown process $X(t)$, was
suggested by Kac \cite{Kac 1956} and firstly studied by McKean
\cite{McKean 1966}. So far numerous works has been done on
McKean-Vlasov type SDEs and applications, see for example, Ahmed
\cite{Ahmed 2007}, Ahmed-Ding \cite{Ahmed-Ding 1995}, Borkar-Kumar
\cite{Borkar-Kumar 2010}, Chan \cite{Chan 1994}, Crisan-Xiong
\cite{Crisan-Xiong 2010}, Kotelenez \cite{Kotelenez 1995},
Potelenez-Kurtz \cite{Kotelenez-Kurtz 2010}, and so on. It is worthy
to point out that (\ref{eq:1.1}) is a particular case of the
following general version,
\begin{eqnarray}\label{eq:1.2}
X(t)&=&x+\int_0^tb(s,X(s),\dbE\phi ^b\mathbf{[}s,X(s),\xi \mathbf{]}%
_{\xi =X(s)})ds \nonumber \\
 \quad\quad &&+\int_0^t\sigma (s,X(s),\dbE\phi ^\sigma
\mathbf{[}s,X(s),\xi \mathbf{]}_{\xi =X(s)})dW_s,
\end{eqnarray}
which can be regarded as a natural generalization of classical SDEs.
Mathematical mean field approaches play a crucial role in diverse
areas, such as physics, chemistry, economics, finance and games
theory, see for example Lasry--Lions \cite{Lasry-Lions 2007}, Dawson
\cite{Dawson 1983}, Huang-Malhame-Caines \cite{Huang-Malhame-Caines
2006}. In a recent work of Buckdahn--Djehiche--Li--Peng
\cite{Buckdahn-Djehiche-Li-Peng 2009}, a notion of mean-field
backward stochastic differential equation (MF-BSDE for short) of the
form
\begin{eqnarray*}\label{eq:1.3}
Y_t=\xi
+\dint_t^T\dbE'f(s,\omega,\omega',Y_s(\omega),Z_s(\omega),Y_s(\omega'),Z_s(\omega'))ds-\dint_t^TZ_sdW_s
,
\end{eqnarray*}
with $t\in[0,T]$ was introduced to investigate one special
mean-field problem in a purely stochastic approach.

In this paper, we would like to introduce mean-field backward doubly
stochastic differential equation (MF-BDSDE for short)
\begin{eqnarray}\label{eq:1.4}
\nonumber Y_t &=&\xi +\int_t^Tf(s,Y_s,Z_s,\Gamma ^f(s,Y_s,Z_s))ds \\
&&+\int_t^Tg(s,Y_s,Z_s,\Gamma ^g(s,Y_s,Z_s))d\overleftarrow{B}_s-\int_t^TZ_sd%
\overrightarrow{W}_s,
\end{eqnarray}
where
\begin{eqnarray}\label{1.4}
[\Gamma^l(s,Y_s,Z_s)](\omega)=\int_\Omega \theta^l(s,\omega ,\omega
^{\prime },Y_s(\omega ),Z_s(\omega ),Y_s(\omega'),Z_s(\omega'
))\dbP(d\omega ^{\prime }),
\end{eqnarray}
with $l=f,$ $g.$ For convenience, we also denote $$\dbE^{\prime
}[\theta ^l(s,Y_s,Z_s,Y_s^{\prime },Z_s^{\prime
})]:=\Gamma^l(s,Y_s,Z_s),$$ when there is no abuse of notation.
Following the basic ideas in \cite{Pardoux-Peng 1994}, we firstly
discuss the existence and uniqueness of solutions for MF-BDSDE
(\ref{eq:1.4}), which obviously extends the results in both
\cite{Pardoux-Peng 1994} and \cite{Buckdahn-Li-Peng 2009}. It is
worthy to point out that MF-BDSDEs not just is a natural
generalization of BSDEs and MF-BSDEs from the view of mathematics.
Our study on them also is motivated by the problems in the following
two aspects.

As is well-known to us, the study on stochastic partial differential
equations have increasingly been a popular issue in recent years. As
one kind of them, stochastic partial differential equations of
McKean-Vlasov type were discussed in \cite{Kotelenez-Kurtz 2010}. In
fact, such equations were obtained as continuum limit from empirical
distribution of a large number of SDEs, coupled with mean-field
interaction. We also refer the reader to \cite{Crisan-Xiong 2010}
and \cite{Kotelenez 1995} for more details along this. On the other
hand, we would also like to mention the work of Buckdahn--Li--Peng
\cite{Buckdahn-Li-Peng 2009} who studied one kind of nonlocal
deterministic PDEs. In virtue of the "backward semigroup" method
they obtained the existence and uniqueness of viscosity solution for
nonlocal PDEs via mean-field BSDEs (\ref{eq:1.3}) in a Markovian
framework and McKean-Vlasov forward equation. Motivated by the above
two cases, in this paper we will give some discussions on one kind
of nonlocal stochastic partial differential equations. Since our
backward equation here is allowed to dependent on
$Z^{0,x_0}(\cdot)$, therefore the nonlocal SPDEs here is not a
direct generalization of deterministic PDEs in
\cite{Buckdahn-Li-Peng 2009} to the stochastic case. Some additional
necessary and essential terms are required in our SPDE to meet the
general case here, see (\ref{eq:4.3}) below. On the other hand,
comparing with the case in \cite{Pardoux-Peng 1994} and
\cite{Buckdahn-Li-Peng 2009}, due to the nonlocal property of SPDEs
and the interesting measurability of the corresponding solution
$u(t,x)$, one fundamental and important term $u'(t,x)$ is required
to meet such general case, see also Remark 4.2 below. Instead of
investigating the limit result for such equations, we will study the
SPDEs from other aspects. A probabilistic interpretation for the
solution to such kind of SPDEs is derived by a connection between
them and decoupled forward-backward doubly differential equations of
mean-field type, which extends the results in \cite{Pardoux-Peng
1994} to the mean-field case.

The second motivation stems from the study of optimal control
problem and certain stochastic differential games problems. Some
related works along this have followed two main venues. On the one
hand, optimal control of mean-field (forward) stochastic
differential equations was discussed in Andersson--Djehiche
\cite{Andersson-Djehiche 2011}, Buckdahn--Djehiche--Li
\cite{Buckdahn-Djehiche-Li 2011} and Meyer-Brandis--Oksandal--Zhou
\cite{Meyer-Brandis-Oksandal-Zhou 2011} where stochastic maximum
principle were derived as a necessary condition of the optimal
control. On the other hand, optimal control problem for backward
doubly stochastic differential system (or forward-backward system)
were spread out in \cite{Bahlali-Gherbal 2010}, \cite{Han-Peng-Wu
2010}, \cite{Zhang-Shi 2010} where the corresponding linear
quadratic problem and nonzero sum stochastic differential games were
also investigated. Inspired by above two case, it is natural for us
to consider the optimal control problem for backward doubly
stochastic system of mean-field type. Since the terms $\Gamma^l$
with $l=f,g$ (see (\ref{eq:5.1})) have more general feature than the
corresponding one in \cite{Andersson-Djehiche 2011},
\cite{Buckdahn-Djehiche-Li 2011} and
\cite{Meyer-Brandis-Oksandal-Zhou 2011}, we have to introduce
$\dbE^*$ in our adjoint equation, which is slight different from
$\dbE'$. Note that such kind of skill also appears in
\cite{Shi-Wang-Yong 2011}. On the other hand, since $\theta^f$,
$\theta^g$ and $l$ are allowed to depend on $u(\cdot,\omega')$, some
new terms appear and is expressed by $\dbE^*$ in our maximum
principle (see (\ref{eq:5.9})) which is totally different from all
the previous literature. We believe that such new features will lead
to some other interesting things and we hope to explore and study
them carefully in the future.

The paper is organized as follows. In Section  \ref{sec:2}, we will
present some preliminary notations needed in the whole paper. In
Section \ref{sec:3}, we consider the existence and uniqueness of
solution for MF-BDSDE (\ref{eq:1.4}). In Section \ref{sec:4}, we
give the probabilistic interpretation for the solutions to a class
of nonlocal SPDEs by means of MF-BDSDE. In Section \ref{sec:5}, we
discuss one optimal control problem of MF-BDSDE. In Section
\ref{sec:6}, we investigate one backward doubly stochastic LQ
problem of mean-field type to show one direct application of result
in Section \ref{sec:5}.

\section{Preliminaries}\label{sec:2}

Throughout this paper, let $\left( \Omega ,\mathcal{F},\dbP\right) $
be a complete probability space on which are defined
two mutually independent Brownian motion $%
\left\{ W_t\right\} _{t\geq 0}$ and $\left\{ B_t\right\} _{t\geq
0}$, with value respectively in $\dbR^d$ and $\dbR^l.$ We denote by
\[
\mathcal{F}_t\doteq {\cal  F}_t^W\vee \mathcal{F}_{t,T}^B,\ \forall
t\in \left[ 0,T\right] ,
\]
where ${\mathcal  N}$ is the class of $\dbP$-null sets of $%
{\cal F}$ and
\begin{eqnarray*}
{\mathcal F}_t^W\doteq \sigma \left\{ W_r;\ 0\leq r\leq t\right\}
\vee {\mathcal  N},\quad {\cal  F}_{t,T}^B\doteq \sigma \left\{
B_T-B_r;\ t\leq r\leq T\right\} \vee {\cal  N}.
\end{eqnarray*}
In this case, the collection $\left\{ {\cal  F}%
_t,t\in \left[ 0,T\right] \right\} $ is neither increasing nor
decreasing, while $\left\{ {\cal  F}_t^W;t\in \left[ 0,T\right]
\right\} $ is an increasing filtration and $\left\{ {\cal
F}_{t,T}^B;t\in \left[ 0,T\right] \right\} $ is a decreasing
filtration.

Let $(\Omega^2,{\cal F}^2,\dbP^2)=(\Omega\times\Omega,{\cal
F}\otimes {\cal F},\dbP\otimes \dbP)$ be the completion of the
product probability space of the above $(\Omega,{\cal F},\dbP)$ with
itself, where we define ${\cal F}^2_t={\cal F}_t\otimes{\cal F}_t$
with $t\in[0,T]$ and ${\cal F}_t\otimes{\cal F}_t$ being the
completion of ${\cal F}_t\times{\cal F}_t$. It is worthy of noting
that any random variable $\xi=\xi(\omega)$ defined on $\Omega$ can
be extended naturally to $\Omega^2$ as
$\xi'(\omega,\omega')=\xi(\omega)$ with
$(\omega,\omega')\in\Omega^2$. For $H=\dbR^n$, ect, let
$L^1(\Omega^2,{\cal F}^2,\dbP^2;H)$ be the set of random variable
$\xi:\Omega^2\rightarrow H$ which is ${\cal F}^2$-measurable such
that
$\dbE^2|\xi|\equiv\int_{\Omega^2}|\xi(\omega',\omega)|\dbP(d\omega')\dbP(d\omega)<\infty.$
For any $\eta\in L^1(\Omega^2,{\cal F}^2,\dbP^2;H),$ we denote
$$\dbE'\eta(\omega,\cdot)=\int_\Omega\eta(\omega,\omega')\dbP(d\omega'),\quad
\dbE^{*}\eta(\cdot,\omega)=\int_\Omega\eta(\omega',\omega)\dbP(d\omega').$$
Particularly, for example, if $\eta_1(\omega,\omega')=\eta(\omega')$
and $\eta_2(\omega,\omega')=\eta_2(\omega)$, then
\begin{eqnarray*}
\dbE'\eta_1=\int_\Omega\eta_1(\omega')\dbP(d\omega')=\dbE\eta_1,\quad
\dbE^*\eta_2=\int_\Omega\eta_2(\omega)\dbP(d\omega)=\dbE\eta_2.
\end{eqnarray*}
Hence, in what follows, $\dbE'$ and $\dbE^*$ will be used when we
need to distinguish $\omega'$ from $\omega$, which is the case of
both $\omega$ and $\omega'$ appearing at the same time. On the one
hand, the well definition of $\dbE'$ above gives the precise meaning
of $\Gamma^f$ and $\Gamma^g$ in (\ref{1.4}). On the other hand, it
also indicates that, for example, the operator $\Gamma^f$ is {\it
nonlocal} in the sense that the value
$\Gamma^f(s,\omega,Y(s,\omega),Z(s,\omega))$ of
$\Gamma^f(s,Y(s),Z(s))$ at $\omega$ depends on the whole set
$$\{(Y(s,\omega'),Z(s,\omega'))\bigm|\omega'\in\Omega\},$$
not just on $(Y(s,\omega),Z(s,\omega))$.

At last, we would like to introduce some spaces of functions
required in the sequel.
\begin{eqnarray*}
S^2\left([0,T];\dbR^n\right)&=&\{\varphi:[0,T]\times\Omega\rightarrow
\dbR^n, \Bigm| \varphi_t \ \mbox{is} \ \mathcal{F}_t \
\mbox{-measurable process}\\
&&
\mbox{such that}\ \dbE(\sup_{0\leq t\leq T}|\varphi_{t}|^{2})<\infty\},\\
M^2(0,T;\dbR^n)&=&\{\varphi_t:[0,T]\times\Omega\rightarrow \dbR^n,
\Bigm| \varphi_t\
 \mbox {is }\ \mathcal{F}_t\mbox{-measurable process}\\
&&
\mbox{such that} \  \dbE\dint_{0}^T|\varphi_{t}|^{2}dt<\infty\},\\
L^{2}(\Omega, {\cal
F}_T,P;\dbR^{n})&=&\{\xi:[0,T]\times\Omega\rightarrow \dbR^n, \Bigm|
\xi\ \mbox{is}\ {\cal
F}_T \mbox{-measurable  random }\\
 && \mbox{variable such that}\ \dbE|\xi|^2<\infty\}.
\end{eqnarray*}
\section{The unique solvability of MF-BDSDEs}\label{sec:3}

In this section, we will discuss the existence and uniqueness of
adapted solution for MF-BDSDE (\ref{eq:1.4}) which is rewritten
below (for convenience):
\begin{eqnarray}\label{eq:3.1}
\nonumber Y_t &=&\xi +\int_t^Tf(s,Y_s,Z_s,\Gamma ^f(s,Y_s,Z_s))ds \\
&&+\int_t^Tg(s,Y_s,Z_s,\Gamma ^g(s,Y_s,Z_s))d\overleftarrow{B}_s-\int_t^TZ_sd%
\overrightarrow{W}_s,
\end{eqnarray}
with $l=f,$ $g.$ Before it, we make the following assumptions.

\begin{enumerate}\item[(H1)]
\begin{enumerate}
\item[(i)]
$\xi \in L^{2}\left( \Omega ,{\cal F}_T,\dbP;\dbR^{n}\right)$. $f:
\Omega \times [0,T] \times \dbR^{n} \times \dbR^{n \times d} \times
\dbR^{k_1} \rightarrow\dbR^{n}$ is measurable and for all
$(t,y,z,\gamma)\in [0,T]\times \dbR^{n+n\times d+k_1}$,
$(t,\omega)\mapsto f(t,\omega, y,z,\gamma)$ is ${\cal
F}_t$-measurable. $g: \Omega \times [0,T] \times \dbR^{n} \times
\dbR^{n \times d} \times \dbR^{k_2}\rightarrow \dbR^{m \times l}$ is
measurable and for all $(t,y,z,\gamma)\in [0,T]\times
\dbR^{n+n\times d+k_2}$, $(t,\omega)\mapsto g(t,\omega, y,z,\gamma)$
is ${\cal F}_t$-measurable. The map
$\theta^l:[0,T]\times\Omega^2\times\dbR^{2n+2n\times2d}\to\dbR^m$ is
measurable and for all $(t,y,z,y'$,
$z')\in[0,T]\times\dbR^{2n+2n\times2d}$, the map
$(t,\omega,\omega')\mapsto\theta^l(t,\omega,\omega',y,z,y',z')$ is
${\cal F}^2_t$-measurable on $[0,T]$, $l=f,g$.
\item[(ii)]$f$ and $g$ satisfy uniformly Lipschitz condition to $(y,z,\gamma)$,
that is, there exist positive constants $L_i$ $K_i$ and $\alpha_j$
with $i=y,z,y',z',\gamma,$ $j=1,2,3,4$, such that
\begin{eqnarray*}
&&|f(t,y_{1},z_{1},\gamma_{1})-f(t,y_{2},z_{2},\gamma_{2})|\\
 &\leq& L_y|y_{1}-y_{2}|+L_z|z_{1}-z_{2}|+L_{\gamma}|\gamma_{1}-\gamma_{2}|,\\
&&|g(t,y_{1},z_{1},\gamma_{1})-g(t,y_{2},z_{2},\gamma_{2})|^2\\
 &\leq&
K_y^2|y_{1}-y_{2}|^2+\alpha_1|z_{1}-z_{2}|^2+\alpha_2|\gamma_{1}-\gamma_{2}|^2,\\
&&|\theta^f(t,\omega,\omega',y_1,z_1,y_1',z_1')-\theta^f(t,\omega,\omega',y_2,z_2,y_2',z_2')|\\
&\leq& L_y|y_1-y_2|+L_z|z_1-z_2|+L_{y'}|y_1'-y_2'|+L_{z'}|z_1'-z_2'|,\\
&&|\theta^g(t,\omega,\omega',y_1,z_1,y_1',z_1')-\theta^g(t,\omega,\omega',y_2,z_2,y_2',z_2')|^2\\
&\leq& K^2_y|y_1-y_2|^2+K^2_{y'}|y_1'-y_2'|^2+\alpha_3|z_1-z_2|^2+\alpha_4|z_1'-z_2'|^2,\\
&&\forall(t,\omega,\omega')\in[0,T]\times\Omega^2,\
(y_i,z_i,y_i',z_i')\in\dbR^{4n+2d},i=1,2,
\end{eqnarray*}
and
\begin{eqnarray*}
&&\dbE\int_0^T|\dbE'\theta^l_0(t,\omega,\omega')|^2ds<\infty,\quad
\dbE\int_0^T|l_0(s,\omega)|^2ds<\infty,\ l=f,g,
\end{eqnarray*}
where
$\theta^l_0(t,\omega,\omega')=\theta^l_0(t,\omega,\omega',0,0,0,0),$
$l_0(s,\omega)=l_0(s,\omega,0,0,0).$ Here we assume that
$\alpha_1+\alpha_2\alpha_3+\alpha_2\alpha_4<1$.
\end{enumerate}
\end{enumerate}

\begin{remark}\label{re:3.1}
Under (H1), since we have
\begin{eqnarray*}
&&|\theta^l(t,\omega,\omega',Y(t,\omega),Z(t,\omega),y,z)|\\
&\leq& L(1+|Y(t,\omega)|+|Z(t,\omega)|+|y|+|z|),\ l=f,g,
\end{eqnarray*}
consequently,
\begin{eqnarray*}
|\Gamma^l(t,Y(t),Z(t))|\leq
L(1+|Y(t)|+|Z(t)|+\dbE|Y(t)|+\dbE|Z(t)|),\ l=f,g.
\end{eqnarray*}%
Likewise, for any
$(Y_1(\cdot),Z_1(\cdot)),(Y_2(\cdot),Z_2(\cdot))\in
S^2([0,T];\dbR^{n})\times M^2([0,T];\dbR^{n\times d})$,
\begin{eqnarray*}
&&|\Gamma^f(t,Y_1(t),Z_1(t))-\Gamma^f(t,Y_2(t),Z_2(t))|\\
&\leq& L_y|Y_1(t)-Y_2(t)|+L_z|Z_1(t)-Z_2(t)|\\
&&+L_{y'}\dbE|Y_1(t)-Y_2(t)|+L_{z'}\dbE|Z_1(t)-Z_2(t)|,\\
&&|\Gamma^g(t,Y_1(t),Z_1(t))-\Gamma^g(t,Y_2(t),Z_2(t))|^2\\
&\leq& K^2_y|Y_1(t)-Y_2(t)|^2+K_{y'}^2\dbE|Y_1(t)-Y_2(t)|^2\\
&&+\alpha_3|Z_1(t)-Z_2(t)|^2+\alpha_4\dbE|Z_1(t)-Z_2(t)|^2.
\end{eqnarray*}
The above two estimates will play an interesting role in the next
theorem.
\end{remark}
\begin{theorem}\label{thm:3.1}
Suppose (H1) holds. Then MF-BDSDE {\rm (\ref{eq:3.1})} admits a
unique adapted solution $(Y,Z)\in S^2\left([0,T];\dbR^n\right)
\times M^2(0,T;\dbR^{n\times d})$.
\end{theorem}

{\bf Proof.} For any $(y,z)\in S^2(0,T;\dbR^n)\times
M^2(0,T;\dbR^{n\times d})$, we consider the following MF--BDSDE
\begin{eqnarray*}
Y_t &=&\xi +\int_t^Tf(s,Y_s,Z_s,\dbE^{\prime }[\theta
^l(s,Y_s,Z_s,y_s^{\prime
},z_s^{\prime })])ds \\
&&\ +\int_t^Tg(s,Y_s,Z_s,\dbE^{\prime }[\theta
^l(s,Y_s,Z_s,y_s^{\prime },z_s^{\prime
})])d\overleftarrow{B}_s-\int_t^TZ_sd\overrightarrow{W}_s.
\end{eqnarray*}
According to Theorem 1.1 in \cite{Pardoux-Peng 1994},
 there exists a unique pair of solution
$(Y,Z)\in S^2\left([0,T];\dbR^n\right)\times M^2(0,T;\dbR^{n\times
d})$. Hence, if we define $\Theta (y,z)=(Y,Z),$ then $\Theta $ maps
from $S^2\left([0,T];\dbR^n\right)\times M^2(0,T;\dbR^{n\times d})$
to itself. We now show that $\Theta $ is contractive. To this end,
take any $(y^i,z^i)\in S^2\left([0,T];\dbR^n\right) \times
M^2(0,T;\dbR^{n\times d})$, $(i=1,2)$, and let
\[
(Y^i,Z^i)=\Theta (y^i,z^i).
\]
We denote by $(\widehat{Y},\widehat{Z})=(Y^1-Y^2$,$Z^1-Z^2)$ and $(\widehat{y%
},\widehat{z})=(y^1-y^2$,$z^1-z^2).$ Using It\^{o} formula to $e^{\beta t}|%
\widehat{Y}_t|^2$ we have
\begin{eqnarray*}
&&\dbE e^{\beta t}|\widehat{Y}_t|^2+\dbE\int_t^Te^{\beta s}|\widehat{Z}%
_s|^2ds+\dbE\int_t^T\beta e^{\beta s}|\widehat{Y}_s|^2ds \\
&=&2\dbE\int_t^Te^{\beta
s}\widehat{Y}_s\widehat{f}(s)ds+\dbE\int_t^Te^{\beta
s}|\widehat{g}(s)|^2ds,
\end{eqnarray*}
where
\begin{eqnarray*}
\widehat{f}(s) &=&f(s,Y^1_s,Z^1_s,\dbE^{\prime }[\theta ^f(s,Y^1_s,Z^1_s,y_s^{1\prime},z_s^{1\prime })])\\
&&-f(s,Y^2_s,Z^2_s,\dbE^{\prime }[\theta^f(s,Y^2_s,Z^2_s,y_s^{2\prime},z_s^{2\prime })]),\\
\widehat{g}(s) &=&g(s,Y^1_s,Z^1_s,\dbE^{\prime }[\theta ^g(s,Y^1_s,Z^1_s,y_s^{1\prime},z_s^{1\prime })])\\
&&-g(s,Y^2_s,Z^2_s,\dbE^{\prime
}[\theta^g(s,Y^2_s,Z^2_s,y_s^{2\prime},z_s^{2\prime })]).
\end{eqnarray*}
Hence from (H1), Remark \ref{re:3.1} above and the inequality $ab
\leq\displaystyle\frac{1}{\delta}a^2+\delta b^2$, we have
\begin{eqnarray*}
&&\dbE e^{\beta t}|\widehat{Y}_t|^2+\dbE\int_t^Te^{\beta s}|\widehat{Z}%
_s|^2ds+\dbE\int_t^T\beta e^{\beta s}|\widehat{Y}_s|^2ds \\
&\leq &2\dbE\int_t^Te^{\beta s}[(1+L_\gamma
)|\widehat{Y}_s|(L_y|\widehat{Y}_s|+L_z|\widehat{Z}_s|)+L_\gamma |\widehat{Y}_s|(L_{y^{\prime }}\dbE|\widehat{y}%
_s|+L_{z^{^{\prime }}}\dbE|\widehat{z}_s|)]ds \\
&&+\dbE\int_t^Te^{\beta s}[K_y^2(1+\alpha
_2)|\widehat{Y}_s|^2+(\alpha
_1+\alpha _2\alpha _3)|\widehat{Z}_s|^2+\alpha _2K_{y^{\prime }}\dbE|\widehat{y}%
_s|^2+\alpha _2\alpha _4\dbE|\widehat{z}_s|^2]ds \\
&\leq &M_1\dbE\int_t^Te^{\beta s}|\widehat{Y}_s|^2ds+M_2\dbE\int_t^Te^{\beta s}|%
\widehat{Z}_s|^2ds\\
&&+M_3\dbE\int_t^Te^{\beta s}|\widehat{y}_s|^2ds+M_4\dbE\int_t^Te^{%
\beta s}|\widehat{z}_s|^2ds,
\end{eqnarray*}
where
\begin{eqnarray*}
M_1 &=&K_y^2(1+\alpha _2)+(1+L_\gamma )(L_y+L_z\cdot
C+\frac{1}{2}L_\gamma ^2+L_\gamma
L_{z^{^{\prime }}}\cdot C), \\
M_2 &=&\frac{(1+L_\gamma )L_z}C+(\alpha _1+\alpha _2\alpha _3), \\
M_3 &=&\frac{1}{2}L_{y^{\prime }}^2+\alpha _2K_{y^{\prime }};\quad
M_4=\frac{L_\gamma L_{z^{^{\prime }}}}C+\alpha _2\alpha _4.
\end{eqnarray*}
After some simple calculations, it is easy to see that
\begin{eqnarray*}
&&\dbE\int_t^T\frac{\beta-M_1}{1-M_2}e^{\beta s}|\widehat{Y}_s|^2ds+\dbE\int_t^Te^{\beta s}|\widehat{Z}%
_s|^2ds\\
&\leq & \frac{M_4}{1-M_2}\left[\dbE\int_t^T\frac{M_3}{M_4}e^{\beta
s}|\widehat{y}_s|^2ds
+\dbE\int_t^Te^{\beta s}|\widehat{z}%
_s|^2ds\right].
\end{eqnarray*}
By the assumption imposed on $\alpha_i$, $\Theta$ is a contraction
on $M^2(0,T;\dbR^{n+n\times d})$, thus there is a unique fixed point
$(Y,Z)\in M^2(0,T;\dbR^{n+n\times d})$ which is the solution of
(\ref{eq:3.1}). Moreover, it is easy to check that $Y(\cdot)\in
S^2(0,T;\dbR^n)$. The proof is complete. \hfill $\Box$

\begin{remark}
Note that our result here can fully cover the corresponding results
in {\rm\cite{Buckdahn-Li-Peng 2009}} and {\rm \cite{Pardoux-Peng
1994}}. In fact, if $K_y=\alpha_1=\alpha_2=0$ or
$L_{\gamma}=\alpha_2=0$, our result degenerates respectively to the
case in {\rm\cite{Buckdahn-Li-Peng 2009}} and {\rm\cite{Pardoux-Peng
1994}}.
\end{remark}
Next we introduce a type of forward doubly stochastic differential
equation with mean field type as follows,
\begin{eqnarray}\label{eq:3.2}
\nonumber P_t &=&\eta +\int_0^tf(s,P_s,Q_s,\Gamma ^f(s,P_s,Q_s))ds \\
&&+\int_0^tg(s,P_s,Q_s,\Gamma
^g(s,P_s,Q_s))d\overrightarrow{W}_s-\int_0^tQ_sd
\overleftarrow{B}_s,
\end{eqnarray}
where $\eta$ is $\mathcal{F}_{0}$-measurable. Note that such kind of
equations appear as adjoint equation in the optimal control problem
below. Here we will transform (\ref{eq:3.2}) into the similar form
of (\ref{eq:3.1}). If we define
\begin{eqnarray*}
&&\tilde{B}_t=B_T-B_{T-t},\quad
\tilde{W}_t=W_T-W_{T-t},\quad\tilde{\mathcal{F}}_t\doteq
\mathcal{F}_{t,T}^{\tilde{W}}\vee\mathcal{F}_{t}^{\tilde{B}},\\
&&\tilde{P}_t=P_{T-t}, \quad \tilde{Q}_t=Q_{T-t},\quad t\in[0,T],
\end{eqnarray*}
then $\tilde{\mathcal{F}}_t=\mathcal{F}_{T-t}$, $\eta$ is
$\tilde{\mathcal{F}}_T$-measurable, $\tilde{P}_t$, $\tilde{Q}_t$ are
$\tilde{\mathcal{F}}_t$-measurable, and
\begin{eqnarray}\label{eq:3.3}
\nonumber \tilde{P}_t &=&\eta +\int_t^Tf(s,\tilde{P}_s,\tilde{Q}_s,\Gamma ^f(s,\tilde{P}_s,\tilde{Q}_s))ds \\
&&+\int_t^Tg(s,\tilde{P}_s,\tilde{Q}_s,\Gamma
^g(s,\tilde{P}_s,\tilde{Q}_s))d\overleftarrow{\tilde{W}}_s-\int_t^T\tilde{Q}_sd
\overrightarrow{\tilde{B}}_s.
\end{eqnarray}
Note that (\ref{eq:3.3}) have the similar form as (\ref{eq:3.1}),
then by Theorem 3.1 there exists a unique pair of
$(\tilde{P},\tilde{Q})$ solving (\ref{eq:3.3}), and we have
\begin{theorem}\label{thm:3.2}
Suppose (H1) holds. Then {\rm (\ref{eq:3.2})} admits a unique
adapted solution $(P,Q)\in S^2\left([0,T];\dbR^n\right) \times
M^2(0,T;\dbR^{n\times d})$.
\end{theorem}

\section{Probabilistic interpretation for a class of nonlocal
SPDEs}\label{sec:4}

The connection between BDSDEs and systems of second-order
quasilinear SPDEs was firstly observed by Pardoux and Peng
\cite{Pardoux-Peng 1994}, where the probabilistic interpretation for
second-order SPDEs of parabolic types was derived. Thereafter, this
subject has attracted a lot of research, such as
\cite{Bally-Matoussi 2001}, \cite{Hu-Ren 2009}, \cite{Ren-Lin-Hu
2009}, \cite{Zhang-Zhao 2007}, \cite{Zhang-Zhao 2010}. This section
can be regarded as a continuation of such a theme. In other words,
we will exploit the above theory of MF-BDSDE in order to provide a
probabilistic formula for the solution of a class of nonlocal SPDE.

Given arbitrary $x_0\in\dbR^m$, $(t,x)\in [0,T]\times \dbR^{m}$
being the initial condition, let us consider the following forward
SDE in $\dbR^m$,
\begin{eqnarray}\label{eq:4.1}
X_s^{t,x}=x+\int_t^s\Gamma^b(r,X^{t,x}_r)dr+\int_t^s\Gamma^{\sigma}(r,X^{t,x}_r)d\overrightarrow{W}_r,\quad
s\geq t,
\end{eqnarray}
and the backward equation
\begin{eqnarray}\label{eq:4.2}
\nonumber Y_s^{t,x}&=&\dbE'[h(X_T^{t,x},(X_T^{0,x_0})')]
+\int_s^T\Gamma
_1^f(r,X_r^{t,x},Y_r^{t,x},Z_r^{t,x})dr-\int_s^TZ^{t,x}_rd\overrightarrow{W}_r\\
&+&\int_s^T\Gamma_1^g(r,X_r^{t,x},Y_r^{t,x},Z_r^{t,x})d\overleftarrow{B}_r,\quad
s\geq t,
\end{eqnarray}
where
\begin{eqnarray*}
&&\Gamma^k(r,X^{t,x}_r)=\int_\Omega k(r,X_r^{t,x}(\omega
),X_r^{0,x_0}(\omega ^{\prime }))\dbP(d\omega ^{\prime })
=:\dbE'[k(r,X_r^{t,x},(X_r^{0,x_0})')],
\\
&&\Gamma _1^l(r,X_r^{t,x},Y_r^{t,x},Z_r^{t,x}) \\
\  &=&\int_\Omega \theta^l(r,X_r^{t,x}(\omega ),Y_r^{t,x}(\omega
),Z_r^{t,x}(\omega ),X_r^{0,x_0}(\omega ^{\prime
}),Y_r^{0,x_0}(\omega
^{\prime }),Z^{0,x_0}_r(\omega'))\dbP(d\omega ^{\prime }) \\
&=:&\dbE'[\theta^l(r,X_r^{t,x},Y_r^{t,x},Z_r^{t,x},(X_r^{0,x_0})',(Y_r^{0,x_0})',(Z^{0,x_0}_r)')],
\end{eqnarray*}
with $k=b,\sigma,$ $l=f,g,$ and
\begin{eqnarray*}
&&b:[0,T]\times \dbR^m\times \dbR^m \to \dbR^m,\ \sigma:[0,T]\times
\dbR^m\times \dbR^m \to  \dbR^{m\times d},\\
 &&\theta^f:[0,T]\times
\dbR^m\times \dbR^n\times
 \dbR^{n\times d}\times \dbR^m\times \dbR^n\times\dbR^{n\times d} \to \dbR^n,\\
&&\theta^g :[0,T]\times \dbR^m\times\dbR^n\times
 \dbR^{n\times d}\times \dbR^m\times \dbR^n\times\dbR^{n\times d}
 \to \dbR^{n\times l},\\
 && h:\dbR^m\times \dbR^m \to \dbR^n.
\end{eqnarray*}
It is known that SDE {\rm(\ref{eq:4.1})} has a unique solution if
coefficients satisfy linear growth and the Lipschitz condition, see
\cite{Buckdahn-Li-Peng 2009}. Similarly under suitable assumptions
of $h$, $\theta^f$ and $\theta^g$, we can also obtain naturally the
wellposedness of {\rm(\ref{eq:4.2})} by means of Theorem 3.1 above.

Suppose $\theta^f(t,x,x',\cdot,\cdot,\cdot)$ and
$\theta^g(t,x,x',\cdot,\cdot,\cdot)$ satisfy the conditions in
Theorem \ref{thm:3.1} uniformly for $t$, $x$  and $x'$, and
$\dbE|h(X^{t,x}_T,(X^{0,x_0}_T)')|^2<\infty$, it follows from
Theorem \ref{thm:3.1} that {\rm(\ref{eq:4.2})} admits a unique
solution $(Y_t,Z_t)\in S^2([0,T];\dbR^{n})\times
M^2(0,T;\dbR^{n\times d})$.

We now relate MF-BDSDE {\rm(\ref{eq:4.2})} to the following nonlocal
SPDE:
\begin{eqnarray}\label{eq:4.3}
\nonumber du(t,x)&=&-\left(\mathcal{L}u(t,x)+\dbE'
\theta^f(t,x,u(t,x),\bigtriangledown u(t,x)\cdot
\dbE'[\sigma (t,x,(X_t^{0,x_0})')],\right.\\
\nonumber
&&\left.(X_t^{0,x_0})',u'(t,(X_t^{0,x_0})'),[\bigtriangledown
u'(t,a)\cdot\dbE'\sigma(t,a,(X^{0,x_0}_t)')]\Big|_{a=X^{0,x_0}_t})\right)dt\\
\nonumber &&-\dbE'\theta^g\left(t,x,u(t,x),\bigtriangledown
u(t,x)\cdot
\dbE'[\sigma(t,x,(X_t^{0,x_0})')],(X_t^{0,x_0})', \right.\\
\nonumber &&\left.u'(t,(X_t^{0,x_0})'),[\bigtriangledown
u'(t,a)\cdot\dbE'\sigma(t,a,(X^{0,x_0}_t)')]
\Big|_{a=X^{0,x_0}_t}\right)d\overleftarrow{B}_t, \\  u(T,x)
&=&\dbE'[h(x,(X_T^{0,x_0})')],
\end{eqnarray}
where $u:\dbR\times \dbR^m\times \Omega \rightarrow \dbR^n$, $\mathcal{L}%
u=(Lu_1,\cdots Lu_n)^T,$ with
\begin{eqnarray*}
Lu_i(t,x)&=&\frac 12tr(\dbE'[\sigma(t,x,(X_t^{0,x_0})')]\dbE'[\sigma(t,x,(X_t^{0,x_0})')]^{T}D^2u(t,x))\\
&&+\bigtriangledown u(t,x)\cdot \dbE'[b(t,x,(X_t^{0,x_0})')],
\end{eqnarray*}
and
\begin{eqnarray}\label{eq:4.3'}
 u'(t,(X^{0,x_0}_t)')=u(t,\omega',X^{0,x_0}_t(\omega')),\quad
 t\in[0,T].
\end{eqnarray}
We can assert that
\begin{theorem}\label{thm:4.1}
Suppose that $b$, $\sigma$, $f$ and $g$ satisfy suitable linear
growth and Lipschitz condition, $h_{xx}(x,(X^{0,x_0}_T))$ exists and
$\dbE\dbE'|h(X^{t,x}_T,(X^{0,x_0}_T))|^2<\infty$. Suppose SPDE
{\rm(\ref{eq:4.3})} has a solution $u(t,x)\in C^{1,2}(\Omega\times
[0,T]\times \dbR^m; \dbR^n)$. Then, for any given $(t,x)$, $u(t,x)$
has the following interpretation
\begin{eqnarray}\label{eq:4.4}
u(t,x)=Y^{t,x}_t,\end{eqnarray} where $Y^{t,x}_t$ is determined by
{\rm(\ref{eq:4.1})} and {\rm(\ref{eq:4.2})}. Moreover, the solution
$u(t,x)$ of {\rm(\ref{eq:4.3})} is unique too.
\end{theorem}

{\bf Proof}\ Applying It\^o's formula to $u(t,X_t)$, we obtain
\begin{eqnarray*}
u(T,X_T^{t,x})-u(t,x) &=&\int_t^T[\frac{\partial u}{\partial r}(r,X_r^{t,x})+%
\mathcal{L}u(r,X_r^{t,x})]dr \\
&&+\int_t^T\bigtriangledown u(r,X_r^{t,x})\cdot \dbE'[\sigma
(r,X_r^{t,x},(X_r^{0,x_0})')]d\overrightarrow{W}_r.
\end{eqnarray*}
Because $u(t,x)$ satisfies SPDE {\rm(\ref{eq:4.3})}, it holds that
\begin{eqnarray*}
&&u(T,X_T^{t,x})-u(t,x)\\
&=&-\int_t^T\dbE'\theta^f(s,X_s^{t,x},u(s,X_s^{t,x}),\bigtriangledown
u(s,X_s^{t,x})\cdot
\dbE'[\sigma (s,X_s^{t,x},(X_s^{0,x_0})')],\\
&&(X_s^{0,x_0})',u'(s,(X_s^{0,x_0})'),[\bigtriangledown
u'(s,a)\cdot\dbE'\sigma(s,a,(X^{0,x_0}_s)')]\Big|_{a=X^{0,x_0}_s})ds\\
&&-\int_t^T\dbE'\theta^g(s,X_s^{t,x},u(s,X_s^{t,x}),\bigtriangledown
u(s,X_s^{t,x})
\cdot\dbE'[\sigma(s,X_s^{t,x},(X_s^{0,x_0})')],\\
&&(X_s^{0,x_0})',u'(s,(X_s^{0,x_0})'),[\bigtriangledown
u'(s,a)\cdot\dbE'\sigma(s,a,(X^{0,x_0}_s)')]\Big|_{a=X^{0,x_0}_s})d\overleftarrow{B}_s\\
&&+\int_t^T\bigtriangledown u(s,X_s^{t,x})\cdot \dbE'[\sigma
(s,X_s^{t,x},(X_s^{0,x_0})')]d\overrightarrow{W}_s.
\end{eqnarray*}
By the uniqueness of solution for (\ref{eq:4.1}) and (\ref{eq:4.2}),
it is easy to check that $(u(s,X_s^{t,x}),\bigtriangledown
u(s,X_s^{t,x})\cdot \dbE'[\sigma (s,X_s^{t,x},(X_s^{0,x_0})')]))$
with $s\in[0,T]$ is a solution of {\rm(\ref{eq:4.2})}. Hence it
follows that
\begin{eqnarray*}
u(t,x)= Y^{t,x}_t,\quad t\in[0,T],x\in\dbR^m.
\end{eqnarray*}
Under the above condition, the solution $Y^{t,x}(\cdot)$ is unique,
then the solution $u(t,x)$ of SPDE (\ref{eq:4.3}) is also unique.
\hfill $\Box$
\begin{remark}\label{rmk:4.1}
{\rm(\ref{eq:4.4})} can be regarded a stochastic Feynman-Kac formula
for SPDE {\rm(\ref{eq:4.3})}, which is a useful tool in the study of
the property for SPDE. For example, from Theorem 4.1 we know that
the solution of SPDE {\rm (\ref{eq:4.3})} must be unique if it
exists.
\end{remark}
\begin{remark}
Note that the introduction of function $u'$ in {\rm(\ref{eq:4.3'})}
coincides with the general setting in our discussion. Actually, on
the one hand, comparing with the SPDEs in {\rm \cite{Pardoux-Peng
1994}}, the term $u'(t,(X^{0,x_0}_t))$ is necessary since our SPDE
here is nonlocal. On the other hand, comparing with the case in {\rm
\cite{Buckdahn-Li-Peng 2009}}, here we replace $u(t,(X^{0,x_0}_t)')$
frequently used there with a slight new term $u'(t,(X^{0,x_0}_t)')$
because of the special measurability of $u(t,x)$. As we know, since
$Y^{t,x}(s,\omega)=u(s,\omega,X^{t,x}(s,\omega)),$ thus
$Y^{0,x_0}(s,\omega')=u(s,\omega',X^{0,x_0}(s,\omega'))=u'(s,(X^{0,x_0}_s)')$.
Obviously, if $\theta^g$ equals to zero, i.e, {\rm(\ref{eq:4.3})}
becomes a deterministic PDE, hence $u'(t,(X^{0,x_0}_t)')$ will
degenerates into $u(t,(X^{0,x_0}_t)')$ which is the case in
{\rm\cite{Buckdahn-Li-Peng 2009}}.
\end{remark}
\begin{remark}
In our framework, {\rm(\ref{eq:4.2})} is allowed to depend on
$Z^{0,x_0}(\cdot)$, hence some more necessary terms representing the
nonlocal property of {\rm (\ref{eq:4.3})} are needed. This is
totally different from the discussion in {\rm \cite{Buckdahn-Li-Peng
2009}}.
\end{remark}

\section{An optimal control problem for MF-BDSDEs}\label{sec:5}

In this section, we would like to consider one optimal control
problem for MF-BDSDEs. As a necessary condition, we will derive one
maximum principle. For the reason of simplicity, we assume
$m=n=d=l=k_1=k_2=1$.

Given a convex subset $U\subset \dbR^k$, for any admissible control
$v\in \mathcal{U}_{ad}$, where
\begin{eqnarray*}
\mathcal{U}_{ad}=\left\{v:[0,T]\times\Omega \to U|v\ \mbox{is}\
\mathcal{F}_t\ \mbox{-measurable},
 \dbE\dint_0^T|v_t|^2{\rm d}t<+\infty \right\},
\end{eqnarray*}
 and $\xi \in L^{2}\left( \Omega
,\mathcal{F}_T,P;\dbR\right)$,
 we consider the following MF-BDSDE:
\begin{eqnarray}\label{eq:5.1}
Y^v_t=\xi +\int_t^T\Gamma ^f(s,Y^v_s,Z^v_s,v_s)ds+\int_t^T\Gamma
^g(s,Y^v_s,Z^v_s,v_s)d\overleftarrow{B}_s-\int_t^TZ^v_sd\overrightarrow{W}_s,\nonumber\\
\end{eqnarray}
where $i=f,$ $g,$%
\begin{eqnarray*}
&&\Gamma ^i(s,Y^v_s,Z^v_s,v_s) \\
&=&\int_\Omega \theta^{i}(s,\omega ,\omega ^{\prime },Y^v_s(\omega
),Z^v_s(\omega ),v_s(\omega ),Y^v_s(\omega ^{\prime }),Z^v_s(\omega
^{\prime }),v_s(\omega^{\prime }))\dbP(d\omega ^{\prime }) ,
\end{eqnarray*}
and
\begin{eqnarray*}
&&\theta^f:\Omega^2\times[0,T]\times \dbR\times \dbR\times
 U\times \dbR\times \dbR \times U\to \dbR,\\
&&\theta^g :\Omega^2\times[0,T]\times \dbR\times \dbR\times
 U\times \dbR\times \dbR \times U\to \dbR.
\end{eqnarray*}
The control problem is to find an admissible control to
 minimize over $\mathcal{U}_{ad}$ the cost function of
\begin{eqnarray}\label{eq:5.2}
J(v(\cdot ))=\dbE\int_0^T\Gamma
^l(s,Y^v_s,Z^v_s,v_s)ds+\dbE[\dbE'h(Y^v_0(\omega),Y^v_0(\omega'))],
\end{eqnarray}
where
\begin{eqnarray*}
&&\ \Gamma ^l(s,Y^v_s,Z^v_s,v_s) \\
\ &=&\int_\Omega l(s,\omega ,\omega ^{\prime },Y^v_s(\omega
),Z^v_s(\omega ),v_s(\omega ),Y^v_s(\omega ^{\prime }),Z^v_s(\omega
^{\prime }),v_s(\omega ^{\prime }))\dbP(d\omega ^{\prime })
\end{eqnarray*}
and $h:\Omega^2\times \dbR\times \dbR \rightarrow \dbR,$
\begin{eqnarray*}
l&:& \Omega^2\times[0,T]\times \dbR\times \dbR\times
 U\times \dbR\times \dbR \times U\to \dbR.
\end{eqnarray*}

The following is our main assumptions on the above mappings in this
section:
\begin{enumerate}\item[(H2)]
\begin{enumerate}
\item[(i)]
$\theta^f$, $\theta^g$, $l$ and $h$ are continuous and continuously
differentiable with respect to $y,$ $y',$ $z,$ $z',$ $v,$ $v',$ and
the derivatives of $l$ and $h$ are allowed to be linear growth.
\item[(ii)]
$\theta^f$ and $\theta^g$ satisfy uniformly Lipschitz condition to
$(y,z,y',z',v,v')$, that is, there exist positive constants $L_i$
$K_i$ and $\alpha_j$ with $i=y,z,y',z',v,v',$ $j=3,4$, such that
\begin{eqnarray*}
&&|\theta^f(t,\omega,\omega',y_1,z_1,y_1',z_1',v_1,v_1')
-\theta^f(t,\omega,\omega',y_2,z_2,y_2',z_2',v_2,v_2')|\\
&&\leq L_y|y_1-y_2|+L_z|z_1-z_2|+L_{y'}|y_1'-y_2'|\\
&&+L_{z'}|z_1'-z_2'|+L_v|v_1-v_2|+L_{v'}|v_1'-v_2'|,\\
&&|\theta^g(t,\omega,\omega',y_1,z_1,y_1',z_1',v_1,v_1')
-\theta^g(t,\omega,\omega',y_2,z_2,y_2',z_2',v_2,v_2')|^2\\
&&\leq K^2_y|y_1-y_2|^2+K^2_{y'}|y_1'-y_2'|^2+K^2_v|v_1-v_2|^2\\
&&+K^2_{v'}|v_1'-v_2'|^2+\alpha_3|z_1-z_2|^2+\alpha_4|z_1'-z_2'|^2,\\
&&\forall(t,\omega,\omega')\in[0,T]\times\Omega^2,\
(y_i,z_i,y_i',z_i',v_i,v_i')\in\dbR^{6},i=1,2,
\end{eqnarray*}
and
\begin{eqnarray*}
&&\dbE\int_0^T|\dbE'\theta^l_0(t,\omega,\omega')|^2ds<\infty,\quad\
l=f,g,
\end{eqnarray*}
where
$\theta^l_0(t,\omega,\omega')=\theta^l_0(t,\omega,\omega',0,0,0,0,0,0),$
and $\alpha_3+\alpha_4<1$.
\end{enumerate}
\end{enumerate}
Under the above hypotheses, for every $v(\cdot)\in
\mathcal{U}_{ad}$, by Theorem \ref{thm:3.1}, (\ref{eq:5.1}) admits a
unique strong solution $(Y^v,Z^v)\in S^2(0,T;\dbR)\times
M^2(0,T;\dbR)$, and the cost functional $J$ is well-defined.

Suppose that $\widehat{u}(\cdot)$ is an optimal control and
$(\widehat{Y}(\cdot), \widehat{Z}(\cdot))$ is the corresponding
optimal trajectory. Let $v(\cdot)$ be such that
$\widehat{u}(\cdot)+v(\cdot) \in \mathcal{U}_{ad}$. Since
$\mathcal{U}_{ad}$ is convex, then for any $0\leq \varepsilon\leq
1$, $u^\varepsilon(\cdot)= \widehat{u}(\cdot)+\varepsilon v(\cdot)$
is also in $\mathcal{U}_{ad}$. From Theorem \ref{thm:3.1}, we know
that state equation (\ref{eq:5.1}) has a unique solution, denoted by
$(Y^\varepsilon(\cdot),Z^\varepsilon(\cdot))$ corresponding to
$u^\varepsilon$. Before the main result we require to prove some
basic results.
\begin{lemma}\label{lem:5.1}
Under assumption (H2), for any $t\in[0,T],$ we have
\[
E|Y_t^\varepsilon -\widehat{Y}_t|^2\leq C\varepsilon ^2,\text{ }%
E\int_t^T|Z_s^\varepsilon -\widehat{Z}_s|^2ds\leq C\varepsilon
^2.\text{ }
\]
\end{lemma}

{\bf Proof.}\ Notice that $Y^\varepsilon_t-\widehat{Y}_t$ satisfies
the following MF-BDSDE:
\begin{eqnarray*}
Y^\varepsilon_t-\widehat{Y}_t&=&\dint_t^T[\Gamma^f(s,Y^\varepsilon_s,Z^\varepsilon_s,
u^\varepsilon_s)-\Gamma^f(s,\widehat{Y}_s,\widehat{Z}_s,\widehat{u}_s)]ds\\
&&+\dint_t^T[\Gamma^g(s,Y^\varepsilon_s,Z^\varepsilon_s,
u^\varepsilon_s)-\Gamma^g(s,\widehat{Y}_s,\widehat{Z}_s,\widehat{u}_s)]d\overleftarrow{B}_s\\
&&-\dint_t^T(Z^\varepsilon_s-\widehat{Z}_s)d\overrightarrow{W}_s.
\end{eqnarray*}

Applying It\^o's formula to $|Y^\varepsilon_t-\widehat{Y}_t|^2$, we
have
\begin{eqnarray*}
&&\dbE\left(|Y^\varepsilon_t-\widehat{Y}_t|^2+\dint_t^T|Z^\varepsilon_s-\widehat{Z}_s|^2ds\right) \\
&=&2\dbE\dint_t^T\left<Y^\varepsilon_s-\widehat{Y}_s,\Gamma^f(s,Y^\varepsilon_s,Z^\varepsilon_s
,u^\varepsilon_s)-\Gamma^f(s,\widehat{Y}_s,\widehat{Z}_s,\widehat{u}_s)\right>ds\\
&&+\dbE\dint_t^T|\Gamma^g(s,Y^\varepsilon_s,Z^\varepsilon_s,
u^\varepsilon_s)-\Gamma^g(s,\widehat{Y}_s,\widehat{Z}_s,\widehat{u}_s)|^2ds.
\end{eqnarray*}
From (H2), we have
\begin{eqnarray*}
&&\dbE|Y^\varepsilon_t-\widehat{Y}_t|^2+\dbE\dint_t^T|Z^\varepsilon_s-\widehat{Z}_s|^2ds\\
&\leq&k_1
\dbE\dint_t^T|Y^\varepsilon_s-\widehat{Y}_s|^2ds+k_2\varepsilon^2\dbE\dint_t^T|v_s|^2ds,
\end{eqnarray*}
where $k_i$ ($i=1,2$) are two constants depending on the parameters
in (H2). By Gronwall's inequality the Burkholder-Davis-Gundy
inequality, we obtain the results. \hfill $\Box$

In the following we make the convention that
\begin{eqnarray*}
\widehat{\alpha}(\cdot ) &=&\alpha (\cdot ,\omega ,\omega ^{\prime
}, \widehat{Y}(\cdot ,\omega ),\widehat{Z}(\cdot ,\omega
),\widehat{u}(\cdot,\omega ),
\widehat{Y}(\cdot ,\omega ^{\prime }),\widehat{Z}(\cdot ,\omega^{\prime }),\widehat{u}(\cdot ,\omega ^{\prime })), \\
\alpha ^\varepsilon (\cdot ) &=&\alpha (\cdot ,\omega ,\omega
^{\prime },Y^\varepsilon (\cdot ,\omega ),Z^\varepsilon (\cdot
,\omega ),u^\varepsilon (\cdot ,\omega ),Y^\varepsilon (\cdot
,\omega ^{\prime }),Z^\varepsilon (\cdot ,\omega ^{\prime
}),u^\varepsilon (\cdot ,\omega ^{\prime }))\text{,}
\end{eqnarray*}
where $\omega ,\omega ^{\prime }\in \Omega ,$ $\alpha
=\theta^f,\theta^g,l.$ We introduce the variational equation as
follows,
\begin{eqnarray}\label{eq:5.3}
\xi_t=\psi _t+\int_t^TF_1(s,\xi _s,\eta _s)ds+\int_t^TG_1(s,\xi
_s,\eta _s)d\overleftarrow{B}_s-\int_t^T\eta
_sd\overrightarrow{W}_s,
\end{eqnarray}
where
\begin{eqnarray*}
F_1(s,\xi _s,\eta _s) &=&\dbE^{\prime }[\widehat{\theta^f}_y(s)\xi
_s +\widehat{\theta^f}_z(s)\eta _s+\widehat{\theta^f}_{y^{\prime
}}(s)\xi'
_s+\widehat{\theta^f}_{z^{\prime }}(s)\eta' _s], \\
G_1(s,\xi _s,\eta _s) &=&\dbE^{\prime }[\widehat{\theta^g}_y(s)\xi
_s + \widehat{\theta^g}_z(s)\eta _s+\widehat{\theta^g}_{y^{\prime
}}(s)\xi'_s +\widehat{\theta^g}_{z^{\prime }}(s)\eta'_s],
\end{eqnarray*}
and
\begin{eqnarray*}
\psi (t) &=&\int_t^T\dbE^{\prime }[\widehat{\theta^f}_v(s)v_s
+\widehat{\theta^f}_{v^{\prime }}(s)v'_s]ds +\int_t^T\dbE^{\prime
}[\widehat{\theta^g}_v(s)v_s+\widehat{\theta^g}_{v^{\prime
}}(s)v'_s]d\overleftarrow{B}_s.
\end{eqnarray*}
Here we denote, for example,
\begin{eqnarray*}
&&\dbE^{\prime }[\widehat{\theta^f}_{y}(s)\xi_s]=\int_\Omega
\widehat{\theta^f}_{y}(s,\omega ,\omega ^{\prime })\xi (s,\omega
)\dbP(d\omega ^{\prime })\\
&&\dbE^{\prime }[\widehat{\theta^f}_{y^{\prime
}}(s)\xi'_s]=\int_\Omega \widehat{\theta^f}_{y^{\prime }}(s,\omega
,\omega ^{\prime })\xi (s,\omega ^{\prime })\dbP(d\omega ^{\prime
}).
\end{eqnarray*}
Under (H2), from Theorem 3.1 there exists a unique $(\xi_t,\eta_t)
\in S^2\left([0,T];\mathbb{R}\right)
 \times M^2(0,T;\mathbb{R})$ satisfying (\ref{eq:5.3}).
\begin{lemma}\label{lem:5.2}
 If we denote by
\begin{eqnarray*}
y_t^\varepsilon =\frac{Y_t^\varepsilon -\widehat{Y}_t}\varepsilon -\xi _t,\text{ }%
z_t^\varepsilon =\frac{Z_t^\varepsilon -\widehat{Z}_t}\varepsilon
-\eta _t.
\end{eqnarray*}
Then we have
\begin{eqnarray}\label{eq:5.4}
\lim\limits_{\varepsilon \rightarrow
0}\sup\limits_{t\in[0,T]}\dbE|y_t^\varepsilon
|^2=0,\quad\lim\limits_{\varepsilon \rightarrow
0}\dbE\int_0^T|z_t^\varepsilon |^2dt=0.
\end{eqnarray}
\end{lemma}
{\bf Proof.}\ First we can express $y^\varepsilon$ and
$z^\varepsilon$ as
\begin{equation*}
\left\{
\begin{array}{lll}
-dy^\varepsilon_t&=&\dbE'[f^\varepsilon_{y}(t)y^\varepsilon_t
+f^\varepsilon_z(t)z^\varepsilon_t+f^\varepsilon_{y'}(t)y^{\prime\varepsilon}_t
+f^\varepsilon_{z'}(t)z^{\prime\varepsilon}_t+f^\varepsilon_1(t)]dt\\
&&+\dbE'[g^\varepsilon_y(t)y^\varepsilon_t
+g^\varepsilon_z(t)z^\varepsilon_t+g^\varepsilon_{y'}(t)y^{\prime\varepsilon}_t
+g^\varepsilon_{z'}(t)z^{\prime\varepsilon}_t+g^\varepsilon_1(t)]d\overleftarrow{B}_t\\
&&-z^\varepsilon_td\overrightarrow{W}_t,\\
y^\varepsilon_T&=&0,
\end{array}\right.
\end{equation*}
where, for example, we denote $\delta=f,g$, $
\bar{Y}_{t,\omega}=\widehat{Y}_{t,\omega}+\lambda(Y^{\varepsilon}_{t,\omega}-\widehat{Y}_{t,\omega}),$
$
\bar{u}_{t,\omega}=\widehat{u}_{t,\omega}+\lambda(u^{\varepsilon}_{t,\omega}-\widehat{u}_{t,\omega})$
$$
\delta_{y}^\varepsilon(\cdot)=\dint_0^1
\theta^\delta_y(\cdot,\bar{Y}_{\cdot,\omega},\bar{Z}_{\cdot,\omega},\bar{u}_{\cdot,\omega},
\bar{Y}_{\cdot,\omega'},\bar{Z}_{\cdot,\omega'},\bar{u}_{\cdot,\omega'})d\lambda,
$$
and
\begin{eqnarray*}
\delta_1^{\varepsilon}(\cdot)&=&[\delta_y^\varepsilon(\cdot)-\widehat{\theta^\delta}_y(\cdot)]\xi_{\cdot,\omega}
+[\delta_z^\varepsilon(\cdot)-\widehat{\theta^\delta}_z(\cdot)]\eta_{\cdot,\omega}
+[\delta_{y'}^\varepsilon(\cdot)-\widehat{\theta^\delta}_{y'}(\cdot)]\xi_{\cdot,\omega'}\\
&&+[\delta_{z'}^\varepsilon(\cdot)-\widehat{\theta^\delta}_{z'}(\cdot)]\eta_{\cdot,\omega'}
+[\delta_{v}^\varepsilon(\cdot)-\widehat{\theta^\delta}_{v}(\cdot)]v_{\cdot,\omega}
+[\delta_{v'}^\varepsilon(\cdot)-\widehat{\theta^\delta}_{v'}(\cdot)]v_{\cdot,\omega'}.
\end{eqnarray*}
Applying It\^o's formula to $\left|y^\varepsilon_t\right| ^2$ on
$\left[ t,T\right]$, we get
\begin{eqnarray*}
&&\dbE\left|y^\varepsilon_t\right| ^2
+\dbE\dint_t^T\left|z^\varepsilon_s\right| ^2ds \\
&=&2\dbE\dint_t^T\left\langle y^\varepsilon_s,
\dbE'[f^\varepsilon_y(s)y^\varepsilon_s
+f^\varepsilon_z(s)z^\varepsilon_s+f^\varepsilon_{y'}(s)y^{\prime\varepsilon}_s
+f^\varepsilon_{z'}(s)z^{\prime\varepsilon}_s+f^\varepsilon_1(s)]\right\rangle ds\\
&&+\dbE\dint_t^T\left|\dbE'[g^\varepsilon_y(t)y^\varepsilon_s
+g^\varepsilon_z(s)z^\varepsilon_s+g^\varepsilon_{y'}(s)y^{\prime\varepsilon}_s
+g^\varepsilon_{z'}(s)z^{\prime\varepsilon}_s+g^\varepsilon_1(s)]\right|
^2ds.
\end{eqnarray*}
From (H2), we have
\begin{eqnarray*}
\dbE|y^\varepsilon_t|^2+\dbE\dint_t^T|z^\varepsilon_s|^2ds \leq
k\dbE\dint_t^T|y^\varepsilon_s|^2ds+C_\varepsilon,
\end{eqnarray*}
where $k$ is some suitable constant, $C_\varepsilon\to 0$ as
$\varepsilon\to 0 $.
 By Grownwall's inequality, we obtain the desired result. \hfill $\Box$

Since $\widehat{u}(\cdot)$ is an optimal control, then
\begin{eqnarray}\label{eq:5.5}
\varepsilon^{-1}[J(u^\varepsilon(\cdot))-J(\widehat{u}(\cdot))]\geq
0.
\end{eqnarray}
From this and Lemma \ref{lem:5.2}, we have the following
\begin{lemma}\label{lem:5.3}
Let assumption {\rm(H2)} hold. Then the following variational
inequality holds:
\begin{eqnarray}\label{eq:5.6}
\nonumber &&\dbE\int_0^T\dbE^{\prime }[\widehat{l}_y(s)\xi
_s+\widehat{l}_z(s)\eta _s+\widehat{l}_{y^{\prime
}}(s)\xi'_s+\widehat{l}_{z^{\prime }}(s)\eta'_s+\widehat{l}_v(s)v_s
+\widehat{l}_{v^{\prime }}(s)v'_s]ds \\
&&+ \dbE
\dbE'[h_y(\widehat{Y}_{0,\omega},\widehat{Y}_{0,\omega'})\xi_{0,\omega}
+h_{y'}(\widehat{Y}_{0,\omega},\widehat{Y}_{0,\omega'})\xi_{0,\omega'}]
\geq 0.\nonumber\\
\end{eqnarray}
Here we denote, for example,
\begin{eqnarray*}
\dbE^{\prime }[\widehat{l}_{y^{\prime }}(s)\xi'_s]=\int_\Omega
\widehat{l}_{y^{\prime }}(s,\omega ,\omega ^{\prime })\xi (s,\omega
^{\prime })\dbP(d\omega ^{\prime }).
\end{eqnarray*}
\end{lemma}

{\bf Proof.}\ From the first result of (\ref{eq:5.4}), we derive
\begin{eqnarray*}
&&\dbE\dbE'\varepsilon^{-1}[h(Y^\varepsilon_{0,\omega},Y^\varepsilon_{0,\omega'})
-h(\widehat{Y}_{0,\omega},\widehat{Y}_{0,\omega'})]\\
&=&\dbE\dbE'\varepsilon^{-1}\dint_0^1h_y(\bar{Y}_{0,\omega},
\bar{Y}_{0,\omega'})(Y^\varepsilon_{0,\omega}-\widehat{Y}_{0,\omega})d\lambda\\
&&+\dbE\dbE'\varepsilon^{-1}\dint_0^1h_{y'}(\bar{Y}_{0,\omega},
\bar{Y}_{0,\omega'})(Y^\varepsilon_{0,\omega'}-\widehat{Y}_{0,\omega'})d\lambda\\
&\to&\dbE\dbE'[h_y(\widehat{Y}_{0,\omega},\widehat{Y}_{0,\omega'})\xi_{0,\omega}
+h_{y'}(\widehat{Y}_{0,\omega},\widehat{Y}_{0,\omega'})\xi_{0,\omega'}],\quad
\varepsilon\to 0,
\end{eqnarray*}
where for example,
$\bar{Y}_{0,\omega}=\widehat{Y}_{0,\omega}+\lambda(Y^\varepsilon_{0,\omega}-\widehat{Y}_{0,\omega}).$
Similarly, we have $\varepsilon\to 0$,
\begin{eqnarray*}
&&\varepsilon^{-1}\left\{\dbE\dint_0^T\dbE'[l^\varepsilon(t)-\widehat{l}(t)]dt\right\}\\
&\to&\dbE\int_0^T\dbE^{\prime }[\widehat{l}_y(s)\xi
_s+\widehat{l}_z(s)\eta _s+\widehat{l}_{y^{\prime
}}(s)\xi'_s+\widehat{l}_{z^{\prime }}(s)\eta'_s+\widehat{l}_u(s)v_s
+\widehat{l}_{v^{\prime }}(s)v'_s]ds.
\end{eqnarray*}
Thus (\ref{eq:5.6}) follows. \hfill $\Box$

Now we consider the adjoint equation:
\begin{eqnarray}\label{eq:5.7}
p_t&=&\dbE'h_y(\widehat{Y}_{0,\omega},\widehat{Y}_{0,\omega'})
+\dbE^*h_{y'}(\widehat{Y}_{0,\omega^*},\widehat{Y}_{0,\omega}) \nonumber \\
&&+\int_0^tF_2(s,,p_s,q_s)ds
+\int_0^tG_2(s,p_s,q_s)d\overrightarrow{W}_s-%
\int_0^tq_sd\overleftarrow{B}_s,
\end{eqnarray}
where
\begin{eqnarray*}
F_2(s,,p_s,q_s) &=&\dbE^{\prime }[\widehat{\theta^f}_y(s)p_s
+\widehat{\theta^g}_y(s)q_s+\widehat{l}_y(s)]\\
&&+\dbE^{*}[\widehat{\theta^f}_{y^{'}}(s)p^*_s+\widehat{\theta^g}_{y^{'
}}(s)q^*_s+\widehat{l}_{y^{' }}(s)], \\
G_2(s,p_s,q_s) &=&\dbE^{\prime }[\widehat{\theta^f}_z(s)p_s
+\widehat{\theta^g}_z(s)q_s+\widehat{l}_z(s) ]\\
&&+\dbE^{*}[\widehat{\theta^f}_{z^{'}}(s)p^*_s+\widehat{\theta^g}_{z^{'
}}(s)q^*_s+\widehat{l}_{z^{'}}(s)].
\end{eqnarray*}
Here we denote by, for example,
\begin{eqnarray*}
&&\dbE^{*}\widehat{l}_{y^{'}}(s)=\int_\Omega
\widehat{l}_{y^{'}}(s,\omega ^{*},\omega)\dbP(d\omega ^{*}),\\
&&\dbE^{*}[\widehat{\theta^f}_{y^{'}}(s)p^{*}_s]=\int_\Omega
\widehat{\theta^f}_{y^{'}}(s,\omega ^{*},\omega )p(s,\omega
^{*})\dbP(d\omega ^{*}).
\end{eqnarray*}
The adjoint equation (\ref{eq:5.7}) is a special form of
(\ref{eq:3.2}) with bounded coefficients. Under (H2), it follows
from Theorem 3.2 that (\ref{eq:5.7}) admits a unique solution
$(p_t,q_t)$.

We define the Hamiltonian function $H:[0,T]\times \dbR\times
\dbR\times \dbR\times \dbR \times
\dbR\times\dbR\times\dbR\times\dbR\to \dbR$ as follows:
\begin{eqnarray}\label{eq:5.8}
\nonumber
&&H(t,y_1,z_1,v_1,y_2,z_2,v_2,p,q) \\
&=&\theta^f(t,\omega,\omega',y_1,z_1,v_1,y_2,z_2,v_2)p
+\theta^g(t,\omega,\omega',y_1,z_1,v_1,y_2,z_2,v_2)q  \nonumber \\
&&+l(t,\omega,\omega',y_1,z_1,v_1,y_2,z_2,v_2).
\end{eqnarray}
From variational inequality (\ref{eq:5.6}), we can state the
stochastic maximum principle of optimal control problem for
MF-BDSDEs.

\begin{theorem}\label{thm:5.1}
{\rm(Stochastic maximum principle)}. Let
$(\widehat{Y}(\cdot),\widehat{Z}(\cdot),\widehat{u}(\cdot))$ be an
optimal triple of the control problem $\{(\ref{eq:5.1}),
(\ref{eq:5.2})\}$. Then $\forall v\in U,\ a.e.\ t\in [0,T], a.s.$
\begin{eqnarray}\label{eq:5.9}
[\dbE'\widehat{H}_v(t,\omega,\omega')+\dbE^*H_{v'}(t,\omega^*,\omega)]\cdot(v-\widehat{u}_t)\geq0,
\end{eqnarray}
where for convenience we denote by
\begin{eqnarray}\label{eq:5.10}
&&\widehat{H}(t,\omega,\omega') \nonumber\\
&=&H(t,\omega,\omega',\widehat{Y}_t(\omega),\widehat{Z}_t(\omega),\widehat{u}_t(\omega),
\widehat{Y}_t(\omega'),\widehat{Z}_t(\omega'),
\widehat{u}_t(\omega'),p_t(\omega),q_t(\omega)).\nonumber \\
\end{eqnarray}
\end{theorem}

{\bf Proof.} Applying It\^o's formula to $\langle \xi_t,p_t\rangle$,
we obtain
\begin{eqnarray*}
-\dbE\xi_0p_0&=&\dbE\int_0^T\dbE^{\prime }[l_y(s)\xi _s+ l_z(s)\eta
_s+l_{y^{\prime }}(s)\xi'_s+l_{z^{\prime }}(s)\eta'_s]ds \\
&&-\dbE\int_0^T\dbE^{\prime }[\widehat{\theta^f}_{v'}(s)v'_sp_s+
\widehat{\theta^g}_{v'}(s)v'_sq_s]ds \\
&&-\dbE\int_0^T\dbE^{\prime }[\widehat{\theta^f}_{v}(s)v_sp_s+
\widehat{\theta^g}_{v}(s)v_sq_s]ds.
\end{eqnarray*}
Then by definition of $\dbE^*$ and variational inequality
(\ref{eq:5.6}) above, we have
\begin{eqnarray*}
&&\dbE\int_0^T\dbE^{\prime }[\widehat{\theta^f}_{v'}(s)v'_sp_s+
\widehat{\theta^g}_{v'}(s)v'_sq_s+ \widehat{l}_{v^{\prime
}}(s)v'_s]ds
\\
&&+\dbE\int_0^T\dbE^{\prime }[\widehat{\theta^f}_{v}(s)v_sp_s+
\widehat{\theta^g}_{v}(s)v_sq_s +\widehat{l}_v(s)v_s]ds \geq0.
\end{eqnarray*}
From the definition of Hamiltonian function in (\ref{eq:5.8}),
\begin{eqnarray*}
&&\dbE\int_0^T\left[\dbE^*H_{v'}(t,\widehat{Y}_t(\omega^*),\widehat{Z}_t(\omega^*),\widehat{u}_t(\omega^*),
\widehat{Y}_t(\omega),\widehat{Z}_t(\omega),
\widehat{u}_t(\omega),p_t(\omega^*),q_t(\omega^*)) \right.\\
 &&\left.+
\dbE'H_v(t,\widehat{Y}_t(\omega),\widehat{Z}_t(\omega),\widehat{u}_t(\omega),
\widehat{Y}_t(\omega'),\widehat{Z}_t(\omega'),
\widehat{u}_t(\omega'),p_t(\omega),q_t(\omega)) \right]\nonumber
\cdot v_tdt\geq0.
\end{eqnarray*}
For $\forall v\in U,$ $F$ be an arbitrary element of the
$\sigma$-algebra $\mathcal{F}_t$, set
\[
\overline{v}\left(s\right) =\left\{
\begin{array}{l}
\widehat{u}_s,\qquad s\in [ 0,t) , \\
v,\qquad s\in [ t,,t+\varepsilon ), \ \omega
\in F, \\
\widehat{u}_s, \qquad s\in[t,,t+\varepsilon ), \ \omega\in
\Omega-F,\\
 \widehat{u}_s,\qquad s\in \left[ t+\varepsilon ,T\right]
,
\end{array}
\right.
\]
we have $\overline{v}(s)\in\mathcal{U}_{ad}$. Since $v_t$ satisfies
$\widehat{u}_t+v_t\in \mathcal{U}_{ad}$, then by taking
$v_t=\overline{v}_t-\widehat{u}_t$, we can rewrite above inequality
as
\begin{eqnarray*}
\dbE \mathbf{1}_F \int_t^{t+\varepsilon
}[\dbE'\widehat{H}_v(s,\omega,\omega')+\dbE^*\widehat{H}_{v'}(s,\omega^*,\omega)]\cdot
(v-\widehat{u}_s)ds\geq0,
\end{eqnarray*}
where $\widehat{H}$ is defined in (\ref{eq:5.10}). Differentiating
with respect to $\varepsilon $ at $\varepsilon =0$ gives
\begin{eqnarray*}
\dbE \mathbf{1}_F
[\dbE'\widehat{H}_v(t,\omega,\omega')+\dbE^*\widehat{H}_{v'}(t,\omega^*,\omega)]\cdot
(v-\widehat{u}_t)\geq0,
\end{eqnarray*}
and (\ref{eq:5.9}) holds naturally. \hfill $\Box$

\section{One mean-field backward LQ problem}\label{sec:6}
In this section, we are dedicated to apply the previous maximum
principle to one backward doubly stochastic LQ problem of mean field
type. In this case, by supposing $h(Y_{0,\omega},Y_{0,\omega'}
)=\frac 12Q^1_0Y_{0,\omega} ^2+\frac 12Q^2_0Y_{0,\omega'}^2$ and
\begin{eqnarray*}
&&f(s,\omega ,\omega ^{\prime },Y_{s,\omega} ,Z_{s,\omega}
,v_{s,\omega} ,Y_{s,\omega
^{\prime }},Z_{s,\omega ^{\prime }},v_{s,\omega ^{\prime }}) \\
&=&A_s^1Y_{s,\omega} +B_s^1Z_{s,\omega} +C_s^1v_{s,\omega}
+A_s^2Y_{s,\omega ^{\prime
}}+B_s^2Z_{s,\omega ^{\prime }}+C_s^2v_{s,\omega ^{\prime }}, \\
&&g(s,\omega ,\omega ^{\prime },Y_{s,\omega} ,Z_{s,\omega}
,v_{s,\omega} ,Y_{s,\omega
^{\prime }},Z_{s,\omega ^{\prime }},v_{s,\omega ^{\prime }}) \\
&=&D_s^1Y_{s,\omega} +E_s^1Z_{s,\omega} +F_s^1v_{s,\omega}
+D_s^2Y_{s,\omega ^{\prime
}}+E_s^2Z_{s,\omega ^{\prime }}+F_s^2v_{s,\omega ^{\prime }}, \\
&&l(s,\omega ,\omega ^{\prime },Y_{s,\omega} ,Z_{s,\omega}
,v_{s,\omega} ,Y_{s,\omega
^{\prime }},Z_{s,\omega ^{\prime }},v_{s,\omega ^{\prime }}) \\
&=&\frac 12[M_s^1Y_{s,\omega} ^2+N_s^1Z_{s,\omega}
^2+R_s^1v_{s,\omega}^2+M_s^2Y_{s,\omega ^{\prime
}}^2+N_s^2Z_{s,\omega ^{\prime }}^2+R_s^2v_{s,\omega ^{\prime }}^2],
\end{eqnarray*}
with, for example, $A^i:[0,T]\times\Omega^2\rightarrow \dbR$ being
bounded, $(s,\omega,\omega')\mapsto A_i(s,\omega,\omega')$ being
$\mathcal{F}^2_s$-measurable (such assumption also hold for the
other coefficients), $M^i,R^i$ being nonnegative, $R^i$ being
positive, we can write
 the state equation and the cost functional as
\begin{eqnarray}\label{eq:6.1}
Y_{t,\omega}^v&=&\xi+\int_t^T\{[\dbE^{\prime
}A_s^1]Y_{s,\omega}^v+[\dbE^{\prime
}B_s^1]Z_{s,\omega}^v+[\dbE^{\prime }C_s^1]v_{s,\omega}+\dbE^{\prime }[A_s^2Y_{s,\omega'}^v]\}ds \nonumber\\
&&+\int_t^T\dbE'[B_s^2Z_{s,\omega'}^v+C_s^2v_{s,\omega'}]ds+\int_t^T\{[\dbE^{\prime
}D_s^1]Y_{s,\omega}^v+[\dbE^{\prime }E_s^1]Z_{s,\omega}^v\}d\overleftarrow{B}_{s,\omega} \nonumber\\
&&+\int_t^T\{[\dbE^{\prime }F_s^1]v_{s,\omega} +\dbE^{\prime
}[D_s^2Y_{s,\omega'}^v+E_s^2Z_{s,\omega'}^v+F_s^2v_{s,\omega'}]\}d\overleftarrow{B}_{s,\omega}-\int_t^TZ_{s,\omega}^vd\overrightarrow{W}_{s,\omega}, \nonumber \\
\end{eqnarray}
and
\begin{eqnarray*}
J(v(\cdot )) &=&\frac 12\dbE\left(\int_0^T\dbE^{\prime
}[M_s^1|Y_{s,\omega}^v|^2+N_s^1|Z_{s,\omega}^v|^2+R_s^1v_{s,\omega}^2]ds+\dbE'[Q^1_0|Y_{0,\omega}^v|^2]\right)\\
&&+\frac 12\dbE\left(\int_0^T\dbE^{\prime
}[M_s^2|Y_{s,\omega'}^v|^2+N_s^2|Z_{s,\omega'}^v|^2+R_s^2|v_{s,\omega'}|^2]ds
+\dbE'[Q^2_0|Y_{0,\omega'}^v|^2]\right).
\end{eqnarray*}
For convenience, we write the random coefficients, for example,
$A^i(s.\omega,\omega')$ as $A^i_s$ in the above and the following
part. The Hamiltonian in such setting becomes
\begin{eqnarray*}
&&\ H(s,\omega ,\omega ^{\prime },y_1,z_1,v_1,y_2,z_2,v_2,p,q) \\
\  &=&[A_s^1y_1 +B_s^1z_1 +C_s^1v_1 +A_s^2y_2+B_s^2z_2+C_s^2v_2]p  \\
&&\ +[D_s^1y_1 +E_s^1z_1 +F_s^1v_1 +D_s^2y_2+E_s^2z_2+F_s^2v_2]q  \\
&&\ +\frac
12[M_s^1y_1^2+N_s^1z_1^2+R_s^1v_1^2+M_s^2y_2^2+N_s^2z_2^2+R_s^2v_2^2].
\end{eqnarray*}
It follows from Theorem \ref{thm:5.1} that
\begin{eqnarray}\label{eq:6.2}
0 &=&\dbE^{\prime }[C_s^1p_{s}+F_s^1q_{s}+R_s^1\widehat{u%
}_{s}]+\dbE^{*}[C_s^2p^*_{s}
+F_s^2q^*_{s}+R_s^2\widehat{u}_{s}], \nonumber \\
\end{eqnarray}
where, for example,
\begin{eqnarray*}
&&\dbE^*[C^1_sp_{s}]=\int_\Omega
C^1(s,\omega,\omega')p(s,\omega)\dbP(d\omega'),\\
&&\dbE^*[R^2_s\widehat{u}_s]=\int_\Omega
R^2(s,\omega^*,\omega)u(s,\omega)\dbP(d\omega^*),\\
&& \dbE^*[C^2_sp^*_{s}]=\int_\Omega
C^2(s,\omega^*,\omega)p(s,\omega^*)\dbP(d\omega^*),
\end{eqnarray*}
and
\begin{eqnarray}\label{eq:6.3}
p_t&=&\dbE'[Q^1_0\widehat{Y}_{0,\omega}]
+\dbE^*[Q^2_0\widehat{Y}_{0,\omega}]\nonumber \\
&&+\int_0^tF_2(s,p_s,q_s)ds+\int_0^tG_2(s,p_s,q_s)d\overrightarrow{W}_s-%
\int_0^tq_sd\overleftarrow{B}_s,
\end{eqnarray}
with
\begin{eqnarray*}
F_2(s,,p_s,q_s) &=&\dbE^{\prime }[A_s^1p_s+D_s^1q_s+M_s^1]
+\dbE^{*}[A_s^2p^*_s+D_s^2q^*_s+M_s^2], \\
G_2(s,p_s,q_s) &=&\dbE^{\prime }[B_s^1p_s+E_s^1q_s+N_s^1]
+\dbE^{*}[B_s^2p^*_s+E_s^2q^*_s+N_s^2].
\end{eqnarray*}
\begin{theorem}\label{thm:6.1}
Suppose there exists $\widehat{u}$ satisfies (\ref{eq:6.2}), where
$(p,q)$ satisfy (\ref{eq:6.3}), then it must be the unique optimal
control of above backward LQ problem.
\end{theorem}

{\bf Proof.}
 First we have
\begin{eqnarray*}
&&J(v)-J(\widehat{u}) \\
&=&\frac 12\dbE\int_0^T\dbE^{\prime
}[M_s^1(|Y_{s,\omega}^v|^2-|\widehat{Y}_{s,\omega}|^2)
+N_s^1(|Z_{s,\omega}^v|^2-|\widehat{Z}_{s,\omega}|^2)]ds \\
&&+\frac 12\dbE\int_0^T\dbE^{\prime
}[R_s^1(|v_{s,\omega}|^2-|\widehat{u}_{s,\omega}|^2)+
M_s^2(|Y_{s,\omega'}^v|^2-|\widehat{Y}_{s,\omega'}|^2)]ds \\
&&+\frac 12\dbE\int_0^T\dbE^{\prime
}[N_s^2(|Z_{s,\omega'}^v|^2-|\widehat{Z}_{s,\omega'}|^2)+
R_s^2(|v_{s,\omega'}|^2-|\widehat{u}_{s,\omega'}|^2)]ds \\
&&+\frac 12
\dbE\dbE'[Q^1_0(|Y^v_{0,\omega}|^2-|\widehat{Y}_{0,\omega}|^2)+
Q^2_0(|Y^v_{0,\omega'}|^2-|\widehat{Y}_{0,\omega'}|^2)]\\
&\geq &\dbE\int_0^T\dbE^{\prime
}[M_s^1\widehat{Y}_{s,\omega}(Y_{s,\omega}^v-\widehat{Y}_{s,\omega})
+N_s^1\widehat{Z}_{s,\omega}(Z_{s,\omega}^v-\widehat{Z}_{s,\omega})]ds \\
&&+\dbE\int_0^T\dbE^{\prime
}[R_s^1\widehat{u}_{s,\omega}(v_{s,\omega}-\widehat{u}_{s,\omega})+
M_s^2\widehat{Y}_{s,\omega'}(Y_{s,\omega'}^v-\widehat{Y}_{s,\omega'})]ds \\
&&+\dbE\int_0^T\dbE^{\prime
}[N_s^2\widehat{Z}_{s,\omega'}(Z_{s,\omega'}^v-\widehat{Z}_{s,\omega'})+
R_s^2\widehat{u}_{s,\omega'}(v_{s,\omega'}-\widehat{u}_{s,\omega'})]ds\\
&&+\dbE\dbE'[Q^1_0\widehat{Y}_{0,\omega}(Y_{0,\omega}^v-\widehat{Y}_{0,\omega})
+Q^2_0\widehat{Y}_{0,\omega'}(Y_{0,\omega'}^v-\widehat{Y}_{0,\omega'})].
\end{eqnarray*}
On the other hand, by using It\^{o} formula to
$p_{s,\omega}(Y^v_{s,\omega}-\widehat{Y}_{s,\omega})$ on $[0,T]$, we
have
\begin{eqnarray*}
&&\dbE\dbE'[Q^1_0\widehat{Y}_{0,\omega}(Y_{0,\omega}^v-\widehat{Y}_{0,\omega})
+Q^2_0\widehat{Y}_{0,\omega'}(Y_{0,\omega'}^v-\widehat{Y}_{0,\omega'})]\\
&=&\dbE\int_0^T\dbE^{\prime }[C_s^1p_{s,\omega}+F_s^1q_{s,\omega}](v_{s,\omega}-\widehat{u}_{s,\omega})ds \\
&&-\dbE\int_0^T\dbE^{^{\prime
}}[M_s^1\widehat{Y}_{s,\omega}(Y_{s,\omega}^v-\widehat{Y}_{s,\omega})
+N_s^1\widehat{Z}_{s,\omega}(Z_{s,\omega}^v-\widehat{Z}_{s,\omega})]ds \\
&&+\dbE\int_0^T\dbE^{\prime
}[C_s^2(v_{s,\omega'}-\widehat{u}_{s,\omega'})p_{s,\omega}+
F_s^2(v_{s,\omega'}-\widehat{u}_{s,\omega'})q_{s,\omega}]ds \\
&&-\dbE\int_0^T\dbE^{^{\prime
}}[M_s^2\widehat{Y}_{s,\omega'}(Y_{s,\omega'}^v-\widehat{Y}_{s,\omega'})
+N_s^2\widehat{Z}_{s,\omega'}(Z_{s,\omega'}^v-\widehat{Z}_{s,\omega'})]ds.
\end{eqnarray*}
Thus we have
\begin{eqnarray*}
&&J(v)-J(\widehat{u}) \\
&\geq &\dbE\int_0^T\dbE^{\prime }[C_s^1p_{s,\omega}
+F_s^1q_{s,\omega}](v_{s,\omega}-\widehat{u}_{s,\omega})ds \\
&&+\dbE\int_0^T\dbE^{\prime }[R_s^1u_{s,\omega}(v_{s,\omega}
-\widehat{u}_{s,\omega})+R_s^2u_{s,\omega'}(v_{s,\omega'}-\widehat{u}_{s,\omega'})]ds \\
&&+\dbE\int_0^T\dbE^{\prime }[C_s^2(v_{s,\omega'}
-\widehat{u}_{s,\omega'})p_{s,\omega}+F_s^2(v_{s,\omega'}-\widehat{u}_{s,\omega'})q_{s,\omega}]ds \\
&=&\dbE\int_0^T\dbE^{\prime }[C_s^1p_{s,\omega}
+F_s^1q_{s,\omega}](v_{s,\omega}-\widehat{u}_{s,\omega})ds \\
&&+\dbE\int_0^T[[\dbE^{\prime
}R_s^1]v_{s,\omega}(v_{s,\omega}-\widehat{u}_{s,\omega})
+\dbE^{*}[R_s^2\widehat{u}_{s}](v_{s,\omega}-\widehat{u}_{s,\omega})]ds \\
&&+\dbE\int_0^T\dbE^{*}[C_s^2p^*_{s}
+F_s^2q^*_{s}](v_{s,\omega}-\widehat{u}_{s,\omega})ds.
\end{eqnarray*}
Thus by (\ref{eq:6.2}) we have $J(v)-J(\widehat{u})\geq 0$, which
means $\widehat{u}$ is an optimal control. Since the method of
proving the result of uniqueness is classical and similar to the
case in \cite{Han-Peng-Wu 2010}, we omit it here.


\begin{thebibliography}{00}

\bibitem{Ahmed 2007} N.~U.~Ahmed, \it Nonlinear diffusion governed by McKean-Vlasov
equation on Hilbert space and optimal control, \sl SIAM J. Control
Optim., \rm 46 (2007), 356--378.

\bibitem{Ahmed-Ding 1995} N.~U.~Ahmed and X.~Ding, \it A semilinear McKean-Vlasov stochastic
evolution equation in Hilbert space, \sl Stoch. Proc. Appl., \rm 60
(1995), 65--85.

\bibitem{Andersson-Djehiche 2011} D.~Andersson and B.~Djehiche, \it A maximum principle for SDEs
of mean-field type, \sl Appl. Math. Optim., \rm DOI
10.1007/s00245-010-9123-8.

\bibitem{Bahlali-Gherbal 2010} S.~Bahlali and B.~Gherbal, \it Optimality conditions of controlled backward
doubly stochastic differential equations, \sl Random Oper. Stoch.
Equ.,\rm  18 (2010), 247--265.

\bibitem{Bally-Matoussi 2001} V.~Bally and A.~Matoussi, \it Weak solutions for SPDEs
and backward doubly stochastic differential equations, \sl J.
Theoret. Probab., \rm 14 (2001), 125--164.

\bibitem{Borkar-Kumar 2010} V.~S.~Borkar and K.~S.~Kumar, \it McKean-Vlasov limit in portfolio
optimization, \sl Stoch. Anal. Appl., \rm 28 (2010), 884--906.


\bibitem{Buckdahn-Djehiche-Li 2011} R.~Buckdahn, B.~Djehiche, and J.~Li, \it A general stochastic
maximum principle for SDEs of mean-field type, \sl Appl. Math.
Optim. \rm DOI 10.1007/s00245-011-9136-y.

\bibitem{Buckdahn-Djehiche-Li-Peng 2009} R.~Buckdahn, B.~Djehiche, J.~Li, and S.~Peng, \it Mean-field
backward stochastic differential equations: a limit approach, \sl
Ann. Probab., \rm 37 (2009), 1524--1565.

\bibitem{Buckdahn-Li-Peng 2009} R.~Buckdahn, J.~Li, and S.~Peng, \it Mean-field
backward stochastic differential equations and related partial
differential equations, \sl Stoch. Proc. Appl., \rm 119, (2009)
3133--3154.

\bibitem{Chan 1994} T.~Chan, \it Dynamics of the McKean-Vlasov equation, \sl Ann. Probab. \rm 22 (1994),
431--441.

\bibitem{Crisan-Xiong 2010} D.~Crisan and J.~Xiong, \sl Approximate
McKean-Vlasov representations for a class of SPDEs, \sl Stochastics,
\rm 82 (2010), 53--68.

\bibitem{Dawson 1983} D.~A.~Dawson, \it Critical dynamics and fluctuations for a mean-field model of
cooperative behavior, \sl J. Statist. Phys., \rm 31 (1983), 29--85.

\bibitem{Han-Peng-Wu 2010} Y.~Han, S.~Peng and Z.~Wu, \it Maximum principle for backward doubly stochastic
control systems with applications, \sl SIAM J. Control Optim., \rm
(2010), 4224--4241, .

\bibitem{Hu-Ren 2009} L. ~Hu and Y. ~Ren,
\it Stochastic PDIEs with nonlinear Neumann boundary conditions and
generalized backward doubly stochastic differential equations driven
by L$\acute{\rm e}$vy processes, \sl J. Comput. Appl. Math., \rm 229
(2009), 230--239.

\bibitem{Huang-Malhame-Caines 2006} M.~Huang, R.~P.~Malham\'e, and
P.~E.~Caines, \it Large population stochastic dynamic games:
closed-loop McKean-Vlasov systems and the Nash certainty equivalence
principle, \sl Comm. Inform. Systems, \rm 6 (2006), 221--252.

\bibitem{Kac 1956} M.~Kac, \it Foundations of kinetic theory, \sl Proc. 3rd Berkeley Sympos.
Math. Statist. Prob. \rm 3 (1956), 171--197.



\bibitem{Kotelenez 1995} P.~M.~Kotelenez, \it A class of quasilinear
stochastic partial differential equations of McKean-Vlasov type with
mass conservation, \sl Prob. Theory Rel. Fields, \rm 102 (1995),
159-188.


\bibitem{Kotelenez-Kurtz 2010} P.~M.~Kotelenez and T.~G.~Kurtz, \it
Macroscopic limit for stochastic partial differential equations of
McKean-Vlasov type, \sl Prob. Theory Rel. Fields, \rm 146 (2010),
189--222.

\bibitem{Lasry-Lions 2007} J.~M.~Lasry and P.~L.~Lions, \it Mean field games, \sl Japan J.
Math., \rm 2 (2007), 229--260.

\bibitem{McKean 1966} H.~P.~McKean, \it A class of Markov processes associated with
nonlinear parabolic equations, \sl  Proc. Natl. Acad. Sci. USA, \rm
56 (1966), 1907--1911.

\bibitem{Meyer-Brandis-Oksandal-Zhou 2011} T.~Meyer-Brandis, B.~Oksendal, and X.~Zhou, \it A
mean-field stochastic maximum principle via Malliavin calculus, \sl
A special issue for Mark Davis' Festschrift, \rm to appear in
Stochastics.

\bibitem{Pardoux-Peng 1994} E.~Pardoux and S.~Peng,  \it Backward doubly
stochastic differential equations and systems of quasilinear
parabolic SPDEs, \sl Probab. Theory Related Fields, \rm 98 (1994),
209--227.

\bibitem{Ren-Lin-Hu 2009} Y.~Ren, A.~Lin and L.~Hu, \it Stochastic PDIEs and
backward doubly stochastic differential equations driven by
L$\acute{\rm e}$vy processes,  \sl J. Comput. Appl. Math., \rm 223
(2009), 901--907.


\bibitem{Shi-Wang-Yong 2011} Y.~Shi, T.~Wang and J.~Yong, \it
Mean-field backward stochastic Volterra integral equations, \sl
arxiv:1104.4725v2[math.PR] 5 Jul 2011.



\bibitem{Zhang-Shi 2010} L.~Zhang, and Y.~Shi, \it Maximum Principle for Forward-Backward Doubly Stochastic Control
Systems and Applications, \sl ESAIM: Cont. Optim. Calc. Vari,
 DOI:10.1051/cocv/2010042.


\bibitem{Zhang-Zhao 2007} Q.~Zhang, and H.~Zhao, \it Stationary solutions of SPDEs and infinite
horizon BDSDEs, \sl J. Funct. Anal., \rm 252(2007), 171--219.

\bibitem{Zhang-Zhao 2010} Q.~Zhang, and H.~Zhao, \it Stationary solutions of SPDEs and infinite
 horizon BDSDEs under non-Lipschitz coefficients, \sl J. Differential Equations, \rm 248 (2010) 953--991.


 \end{thebibliography}
\end{document}